\documentclass[12pt]{article}
\usepackage[leqno]{amsmath}
\usepackage{amsthm}
\usepackage{amssymb,amscd,latexsym,stmaryrd}
\usepackage{paralist}
\usepackage[colorlinks,pdfpagelabels,pdfstartview = FitH,bookmarksopen = true,bookmarksnumbered = true,linkcolor = black,plainpages = false,hypertexnames = false,citecolor = black]{hyperref}
\usepackage{bm}
\usepackage{fonttable}                                                         
       
\usepackage{tocloft}
\usepackage{aurical}
\usepackage[T1]{fontenc}
\usepackage[bbgreekl]{mathbbol}
\usepackage[all]{xy}

\newcommand{\C}{{\mathbb{C}}}
\newcommand{\F}{{\mathbb{F}}}

\newcommand{\Pa}{{\mathbb{P}}}
\newcommand{\Q}{{\mathbb{Q}}}
\newcommand{\oQ}{\overline{\Q}}

\newcommand{\R}{{\mathbb{R}}}
\newcommand{\Z}{{\mathbb{Z}}}

\newcommand{\hZ}{\hat{\Z}}

\newcommand{\abb}{\mathrm{ab}}

\newcommand{\car}{\mathrm{char}\,}

\newcommand{\colim}{\mathrm{colim}}

\newcommand{\et}{\mathrm{\acute{e}t}}

\newcommand{\id}{\mathrm{id}}

\newcommand{\inn}{\mathrm{in}}

\newcommand{\rat}{\mathrm{rat}}
\newcommand{\imm}{\mathrm{im}\,}

\newcommand{\oP}{\overline{P}}

\renewcommand{\mod}{\;\mathrm{mod}\;}

\newcommand{\Mor}{\mathrm{Mor}}
\newcommand{\ord}{\mathrm{ord}}

\newcommand{\Quot}{\mathrm{Quot}}

\newcommand{\spec}{\mathrm{spec}\,}

\newcommand{\Aut}{\mathrm{Aut}}

\newcommand{\Gal}{\mathrm{Gal}}

\newcommand{\Hom}{\mathrm{Hom}}

\newcommand{\tP}{\tilde{P}}

\newcommand{\tV}{\tilde{V}}
\newcommand{\tW}{\tilde{W}}

\newcommand{\tvarphi}{\tilde{\varphi}}
\newcommand{\tpsi}{\tilde{\psi}}

\newcommand{\hpi}{\hat{\pi}}

\newcommand{\RRe}{\mathrm{Re}\,}

\renewcommand{\top}{\mathrm{top}}

\newcommand{\Tr}{\mathrm{Tr}}

\newcommand{\tors}{\mathrm{tors}}

\newcommand{\Dh}{{\mathcal D}}
\newcommand{\Eh}{{\mathcal E}}
\newcommand{\Fh}{{\mathcal F}}

\newcommand{\Nh}{\mathcal{N}}
\newcommand{\Oh}{{\mathcal O}}

\newcommand{\emm}{{\mathfrak{m}}}

\newcommand{\eo}{\mathfrak{o}}

\newcommand{\ep}{\mathfrak{p}}

\newcommand{\eX}{\mathfrak{X}}
\newcommand{\ex}{\text{\normalfont\Fontauri x}}

\newcommand{\ey}{\text{\normalfont\Fontauri y}}

\newcommand{\onF}{\overline{F}}
\newcommand{\oK}{\overline{K}}

\newcommand{\oX}{\overline{X}}

\newcommand{\oF}{\overline{\mathbb{F}}}

\newcommand{\ceX}{\check{\eX}}

\newcommand{\deX}{\overset{_{\,\hullet}}{\eX}}

\newcommand{\tf}{\tilde{f}}

\newcommand{\silo}{\xrightarrow{\sim}}

\newcommand{\ent}{\;\widehat{=}\;}
\newcommand{\hullet}{{\scriptscriptstyle \bullet\,}}
\newcommand{\verk}{\mbox{\scriptsize $\,\circ\,$}}

\newtheorem{theorem}{Theorem}[section]
\begin{document}
\title{Primes, knots and periodic orbits}
\author{Christopher Deninger\footnote{Funded by the Deutsche Forschungsgemeinschaft (DFG, German Research Foundation) under Germany's Excellence Strategy EXC 2044--390685587, Mathematics M\"unster: Dynamics--Geometry--Structure and the CRC 1442 Geometry: Deformations and Rigidity} }
\date{\ }
\maketitle
\section{Introduction} \label{sec:1}

We begin by sketching some analogies between number theory and knot theory originally pointed out by Manin, Mazur \cite{Mazur2} and Mumford. Those were later developed by Kapranov, Morishita, Reznikov, Sikora and many other reseachers. We then explain a more structured analogy where the knots arise from the closed orbits of an $\R$-dynamical system. Finally, we explain our recent construction of dynamical systems for arithmetic schemes which realize some aspects - but not all - of these analogies. At least in the beginning we assume more knowledge of analysis than of algebra on the part of the reader.

I am very grateful to Umberto Zannier for the invitation to Pisa and to the organizers of the Colloquium de Giorgi for giving me the occasion to contribute to this series.
\section{Primes and Knots} \label{sec:2}
We first explain how a prime number $p$ can be viewed as a $1$-dimensional object and the ring of integers $\Z$ as $3$-dimensional. For this we need to introduce schemes and the \'etale (Grothendieck) topology.

For a commutative unital ring $R$, the spectrum $\spec R$ consists of the prime ideals $\ep$ of $R$. It is equipped with the Zariski topology whose closed subsets are of the form
\[
V (S) = \{ \ep \in \spec R \mid \ep \supset S \}
\]
with $S$ any subset of $R$. The Zariski topology is almost never Hausdorff and the topological space $\spec R$ is only a weak invariant of $R$. Grothendieck equipped $X = \spec R$ with a sheaf of rings $\Oh_X$. An element $f \in R$ gives a global section of $\Oh_X$ and we may view $f$ as a function on $\spec R$ by sending $\ep$ to the residue class $f \mod \ep$ in $\kappa (\ep) = \Quot (R / \ep)$. Thus we do not have only one field where our functions can take values but many. For example $f \in \Z$ takes values $f ((p)) \in \F_p$ for the prime numbers $p$ and $f ((0)) \in \Q$. In this sense the important idea that ``numbers are functions'' has been realized. However, working with $\spec \Z$ in isolation, this elementary geometric point of view does not give new arithmetic insights. Using functoriality the situation can be improved. Every homomorphism of rings $\varphi : R_1 \to R_2$ induces a continuous map 
\[
\spec R_2 \to \spec R_1\quad \text{via} \; \ep \mapsto \varphi^{-1} (\ep) \; .
\]
One even obtains a map of locally ringed spaces between the ``affine schemes'' $\spec R_i$. Affine schemes can be glued  to give new ``spaces'', the \textit{schemes} of Grothendieck. Moreover, using the entire web of schemes which map to $\spec R$ one can refine the Zariski topologies to Grothendieck-topologies from which very interesting homotopical and cohomological invariants of the original ring $R$ can be obtained. Let us motivate the basic idea behind the \'etale topology:  The usual topology underlying our geometric intuition on a complex manifold has the following property: Any analytic map whose Jacobian is everywhere invertible is itself locally invertible. In the context of algebraic geometry this statement is no longer true for the Zariski topology. The map $x \mapsto x^2$ has the square root as a local inverse around $x = 1$ but the square root is not a polynomial map. One can translate the invertibility of the Jacobian into commutative algebra terms and one arrives at the class of ``etale'' maps of schemes. Without giving the formal definition, let us just note that the map $\spec L \to \spec K$ corresponding to an inclusion of fields $K \subset L$ is \'etale if and only if the extension $L / K$ is finite and separable. For a finite extension of number fields $L / K$ with rings of integers $\eo_L$ and $\eo_K$ the map $\spec \eo_L \to \spec \eo_K$ is \'etale if and only if the extension $L / K$ is unramified at all primes. In general it is \'etale away from the ramified primes.

The above map $x \mapsto x^2$ gives a ring homomorphism $\C [x] \to \C [x]$ and the induced map $\spec \C [x] \leftarrow \spec \C [x]$ is \'etale around the point $\ep = (x-1)$ of $\spec \C [x]$. For a topology closer to our intuition than the Zariski topology, all \'etale maps and in particular the ones in the example should have local inverses. This can only be achieved by generalizing the notion of a topology on a set $X$. In the usual context, the open sets $U$ are subsets of $X$. Instead Grothendieck considered the class of \'etale maps $f : U \to X$ as the ``open sets'' of the \'etale ``topology'' on $X$. Finite intersections are defined by iterated fibre products over $X$, coverings are also easy to define, and many notions of topology can be extended to this more general setting, including cohomology and fundamental groups. Now we have achieved our goal: In the \'etale topology, \'etale maps $f$ are locally invertible! The reason is trivial: Consider the diagram:
\[
\xymatrix{
 & U \ar@{=}[dl]_{\id} \ar[d]^f \\
U \ar[r]^f & X
}
\]
Here the vertical map $f$ defines an ``open set'' on which the horizontal \'etale map $f$ which we want to invert has an inverse, namely the identity on $U$! Historically, in \cite{Serre} Serre had introduced the condition of being ``locally isotrivial'' for fibre bundles on algebraic varieties. This is a weaker condition than local triviality for the Zariski topology. Inspired by this work, Grothendieck invented the \'etale topology and used it to define \'etale cohomology and the \'etale fundamental group.

Let us now look at the inclusion
\[
\spec \F_p \hookrightarrow \spec \Z \; , \; (0) \mapsto (p)
\]
coming from the projection $\Z \to \Z / p = \F_p$. The \'etale fundamental group $\hpi_1 (\spec \F_p)$ is the automorphism group of the universal (pro-)\'etale covering of $\spec \F_p$ i.e. of $\spec \oF_p$. Hence we have
\[
\hpi_1 (\spec \F_p) = \Aut (\oF_p) \; .
\]
This is a (pro-)cyclic group generated by the Frobenius automorphism $x \mapsto x^p$. We may therefore view
\begin{equation}
\label{eq:1}
\hpi_1 (\spec \F_p) = \hZ := \varprojlim_{N} \Z / N
\end{equation}
as an analogue of
\begin{equation}
\label{eq:2}
\pi_1 (S^1) = \Z
\end{equation}
where $S^1$ is the cicle. 

By the way, the ultimate reason why \'etale fundamental groups are always pro-finite like $\hZ = \varprojlim_N \Z / N$ is this: polynomials in one variable have only finite many zeroes. Note that the group $\Z$ in \eqref{eq:2} occurs as the automorphism group of the covering $\R \to S^1 , t \mapsto \exp 2 \pi i t$, and the deck transformations $t \mapsto t + n$ correspond to the infinitely many zeroes $n \in \Z$ of the function $1 - \exp (2 \pi it)$. 

It turns out that the higher \'etale homotopy groups of $\spec \F_p$ vanish so that for the \'etale topology, $\spec \F_p$ is a $K (\hZ , 1)$-space just as $S^1$ is a $K (\Z , 1)$, space for the ordinary topology. Let us now turn to $\spec \Z$. Its Zariski dimension which intuitively corresponds to the complex dimension is one. So its \'etale cohomological dimension which intuitively corresponds to the usual real dimension should be $2$. However the residue fields $\F_p$ of the closed points $(p)$ are not separably closed and the \'etale topological dimension $1$ of the spec $\F_p$'s adds to the \'etale topological dimension of $\spec \Z$. These heuristics suggest that the \'etale cohomogical dimension of $\spec \Z$ should be $2 \cdot 1 + 1 = 3$. A proof using class field theory can be found in \cite{Mazur1}. 
From an arithmetic point of view $\spec \Z$ is not compact, because it only contains the (equivalence classes of the) non-archimedean absolute values of $\Z$, which correspond to the prime ideals $(p)$, but not the archimedean absolute value. The latter is usually denoted by the symbol $\infty$ and we may view $\overline{\spec \Z} = \spec \Z \cup\{ \infty \}$ as a compactification of $\spec \Z$. 

Artin and Verdier extended the \'etale topology in a natural way to $\overline{\spec \Z}$. The resulting cohomology groups differ from the ones of $\spec \Z$ only by $2$-torsion, which for the purposes of this section is harmless. In the next section, the role of $\infty$ will be more important. Intuitively we may view $\overline{\spec \Z}$ as a compact $3$-manifold $M^3$ (without boundary). By a theorem of Minkowski, there are no non-trivial everywhere unramified extensions of the number field $\Q$. In terms of the \'etale topology this result can be restated as $\hpi_1 (\spec \Z) = 1$ which implies $\hpi_1 (\overline{\spec \Z}) = 1$. In this regard $\overline{\spec \Z}$ resembles a closed \textit{simply connected} $3$-manifold, i.e. $S^3$ by the Poincar\'e conjecture. Let us write $K_p$ for the embedding $\spec \F_p \hookrightarrow \overline{\spec \Z}$. Since the image is $(p)$, we may view $K_p$ as $(p)$ embedded into $\overline{\spec \Z}$. As explained above this is reminiscent of an embedding $S^1 \hookrightarrow S^3$ i.e. of a knot $K$. Apparently, this observation was first made by Mumford \cite[Introduction]{Mazur2}.

Here are some analogies between the topological and the arithmetical situation where for simplicity we deal with $\spec \Z$ instead of $\overline{\spec \Z}$. The abelianized fundamental group of the knot complement is cyclic, $\pi^{\abb}_1 (S^3 \setminus K) \cong \Z$. On the other hand, the abelianized \'etale fundamental group of $\spec \Z \setminus K_p$ is the Galois group of the maximal abelian extension of $\Q$ which is unramified away from the prime $p$. By a theorem of Kronecker the latter is the field $\Q (\mu_{p^{\infty}})$ obtained by adjoining all $p^n$-th roots of unity to $\Q$ for $n \ge 1$. Hence we have 
\[
\hpi^{\abb}_1 (\spec \Z \setminus K_p) = \Gal (\Q (\mu_{p^{\infty}}) / \Q) = \hZ^{\times}_p \; .
\]
Here $\hZ_p$ is the ring of $p$-adic numbers. The group $\hZ^{\times}_p$ is not quite (pro $p$) cyclic, but it almost is. The Alexander polynomial of the knot $K$ is the characteristic polynomial of a generator of $\pi^{\abb}_1 (S^3 \setminus K) \cong \Z$ acting on $\pi^{\abb}_1 (\widetilde{S^3 \setminus K}) \otimes_{\Z} \Q$. Here $\widetilde{S^3 \setminus K}$ is the $\Z$-covering of $S^3 \setminus K$ corresponding to the quotient
\[
\pi_1 (S^3\setminus K) \longrightarrow \pi^{\abb}_1 (S^3 \setminus K) \cong\Z \; .
\]
Correspondingly, the quotient
\[
\hpi_1 (\spec \Z \setminus K_p) \longrightarrow \hpi^{\abb}_1 (\spec \Z \setminus K_p) \cong \hZ^{\times}_p
\]
corresponds to a pro-\'etale $\hZ^{\times}_p$-covering $\widetilde{\spec \Z \setminus K_p}$ of $\spec \Z \setminus K_p$. The action of $\hpi^{\abb}_1 (\spec \Z \setminus K_p) \cong \Z^{\times}_p$ on $\hpi^{\abb}_1(\widetilde{\spec \Z \setminus K_p})$ which is the Galois group of the maximal abelian outside of $p$ unramified extension of $\Q (\mu_{p^{\infty}})$, gives rise to the Iwasawa zeta function. Iwasawa thought in terms of Galois groups, the analogy with knot theory via the \'etale topology was observed only later by Mazur \cite{Mazur2}. The work of Iwasawa gave rise to a most fruitful direction in algebraic number theory. Morishita and his co-authors have also applied techniques from knot theory to the study of non-abelian Iwasawa theory \cite{Mo2}. 

Here is another well known analogy. The $\mod 2$ linking number $l (K,L)$ of two knots $K, L$ may be defined by counting $\mod 2$ the number of times that $L$ intersects a Seifert surface with boundary $K$. With this definition it is a nontrivial fact that we have
\begin{equation}
\label{eq:3}
l (K,L) = l (L,K) \quad \text{in} \; \Z / 2 \; .
\end{equation}
For odd prime numbers $p \neq l$ the quadratic residue symbol $\big( \frac{p}{l} \big)$ is $1$ if $p$ is a square $\mod l$ and $-1$ if it is not. If $p$ or $l$ is congruent to $1 \mod 4$ the famous Gauss reciprocity law asserts that
\begin{equation}
\label{eq:4}
\big( \frac{p}{l} \big) = \big( \frac{l}{p} \big) \; .
\end{equation}
Viewing primes as knots, there are good reasons to view $\big( \frac{p}{l} \big)$ as an analogue of $(-1)^{l (K,L)}$. In fact there is a proof of \eqref{eq:3} which can be translated to the arithmetic context via our dictionary and which then yields a proof of \eqref{eq:4}. The limitations of viewing $\overline{\spec\Z}$ as analogous to $S^3$ become apparent though. If both $p$ and $l$ are not congruent to $1 \mod 4$, quadratic reciprocity assert that
\[
\big( \frac{p}{l} \big) = - \big( \frac{l}{p} \big) \; .
\]
According to Kapranov and Smirnov \cite[\S\,3]{KS}, the reason for this assymmetry is the following: the primes $p$ which are not congruent to $1 \mod 4$ should correspond to knots which are not homologous to zero in the $3$-manifold $M^3$ corresponding to $\overline{\spec \Z}$ in our analogy. Taking $S^3$ for $M^3$ is an oversimplification in this regard because in $S^3$ all knots are homologous to zero. In \cite{KS} it is explained how to restore the analogy between quadratic reciprocity and $\mod 2$ linking of knots if more general $3$-manifolds $M^3$ are allowed.

For three knots in $S^3$ which are pairwise unlinked their union can still be a non-trivial link, as for the Borromean rings. This can be shown by the non-triviality of certain triple Massey products for example. In 1938 Redei published a symbol $(p , l, q)$ defined for three prime numbers $p, l, q$ congruent to $1 \mod 4$ whose quadratic residue symbols satisfy $( \frac{p}{l}) = 1$ and $(\frac{p}{q}) = 1$. Morishita observed that non-triviality of the Redei symbol could be interpreted as the three primes being non-trivially linked while being pairwise unlinked. For example, we may view the primes $5, 41, 61$ as an analogue of the Borromean link. He also gave an entirely parallel construction of ``higher Milnor invariants'' both for $n$-tuples of knots and of primes, which reduce to the Redei symbol for $n = 3$, \cite{Mo3}. The relation of Redei symbols with Massey products in Galois cohomology is explained in \cite{Mo4}. Recently \cite{AC} gives an interesting application of Massey products to infinite class field towers. There are very many further analogies between knot theory and prime numbers and more generally between three dimensional topology and algebraic number theory. We refer to the book \cite{Mo1} for a great overview of this branch of mathematics, called arithmetic topology. 
\section{Primes and periodic orbits}\label{sec:3}
In this section we explain certain analogies between number theory and the theory of $\R$-dynamical systems on $3$-manifolds $M^3$ with a $1$-codimensional foliation $\Fh$. Up to isotopy the periodic orbits give embedded circles i.e. knots in $M^3$. The analogies in the present section enhance the ones of section \ref{sec:2}. They are mostly of an analytic nature and relate to analytic number theory contrary to the topological analogies of the previous section. I arrived at these analogies by a long detour via cohomological considerations \cite{D1}. However they can be much more easily motived by the following argument, \cite{Ko1}. For a rational number $f \in \Q^{\times}$ the $p$-adic absolute values of $f$ are
\[
|f |_p = p^{-\ord_p f} \quad \text{if} \; f = p^{\ord_p f} \frac{a}{b} \quad \text{with} \; 0 \neq b , a \in \Z \; \text{prime to} \; p \; .
\]
We also have the ordinary archimedean absolute value $|f|_{\infty} := |f|$. The product formula which immediately follows from these definitions asserts that
\[
\prod_{p \le \infty} |f|_p = 1\; .
\]
Taking the $\log$, we get
\begin{equation}
\label{eq:5}
\sum_{p \neq \infty} \ord_p f \cdot \log p - \log |f|_{\infty} = 0 \; .
\end{equation}
This reminds of, but is more complicated than the formula
\begin{equation}
\label{eq:6}
\sum_{x \in X} \ord_x f = 0 
\end{equation}
for a meromorphic function $f$ on a Riemann surface $X$. In \eqref{eq:6} only integers are added, whereas in \eqref{eq:5} the integers are multiplied with the $\Q$-linearly independent transcendental numbers $\log p$, and there is also the term $\log |f|_{\infty}$ which seems of a different nature. A formula like \eqref{eq:5} but without a term corresponding to $-\log |f|_{\infty}$ was observed by Ghys in the context of Riemann surface laminations. Kopei later showed how to include an analogue of $- \log |f|_{\infty}$: As in \cite{D2} consider triples $(\oX , \Fh, \phi)$ where $\oX$ is a closed smooth $3$-manifold with a $1$-codimensional foliation $\Fh$ by Riemann surfaces, and $\phi$ is a flow, such that each $\phi^t$ maps leaves of $\Fh$ to (possibly other) leaves. The flow should be non-degenerate which implies that there are only finitely many fixed points $x$ and for any $C > 0$ there are only finitely many closed orbits $\gamma$ with length $l (\gamma) \le C$. The fixed points of $\phi$ should lie in finitely many compact leaves. All other leaves should be non-compact and the flow should be transversal to them. Let $X$ be the open complement in $\oX$ of the union over the finitely many compact leaves of $\Fh$. Since the compact leaves contain fixed points, they are invariant under the flow and become sub-dynamical systems of codimension one. The open manifold $X$ is invariant under the flow and foliated by the non-compact leaves of $\Fh$. 

The triple $(X , \Fh , \phi^t)$ can be described more explicitely as follows. Pick a leaf $F$ of $\Fh$ (they are all diffeomorphic), let $Y_{\phi}$ be the vector field generated by the flow and let $\omega_{\phi}$ be the $1$-form on $X$ which is zero on the tangent bundle $T\Fh$ to the foliation and such that $\langle \omega_{\phi} , Y_{\phi} \rangle = 1$. For the period group
\[
\Lambda = \big\{ \int_{\gamma} \omega_{\phi} \mid \gamma \in \pi^{\abb}_1 (X) \big\} \subset \R
\]
we have $\Lambda = \{ t \in \R \mid \phi^t (F) = F \}$. For example, if $\gamma$ is a periodic orbit then $l (\gamma) \in \Lambda$ and $\phi^{l (\gamma)}$ is the Poincar\'e return map on $F$. The group $\Lambda$ acts diagonally on $F \times \R$. Moreover the action is fixed point free and properly discontinuous. The map $F \times \R \to X , (x,t) \mapsto \phi^t (x)$ induces a diffeomorphism of $F \times_{\Lambda} \R$ with $X$. Here the leaves of $\Fh$ correspond to the images of $F \times \{ s \}$ in $F \times_{\Lambda} \R$ and $\phi^t$ becomes translation by $t$ via the second factor. In \cite{ALK1} and \cite{MKNT} the possible shapes of the triples $(\oX , \Fh , \phi^t)$ are determined. For example the Reeb foliation $\Fh$ on $\oX = S^3$ with suitable flows are possibilities. Recall that the Reeb foliation of $S^3$ is obtained by glueing two solid tori $S^1 \times D$ foliated by concentric infinitely extended paraboloids, along their common boundary. Thus there is one compact leaf, the $2$-torus $S^1 \times S^1$ and all the other leaves are diffeomorphic to the plane.

Any choice of a smooth metric on $T \Fh$ over $\oX$ defines conformal structures on the leaves which vary smoothly in the transversal direction. Thus for a smooth $3$-manifold with a smooth $1$-codimensional foliation one obtains the structure of a \textit{Riemann surface lamination} \cite{G}. Let $f$ be a smooth $\Pa^1 (\C)$-valued function whose restrictions to all leaves $F$ of $\Fh$ are holomorphic and let $d_{\Fh} f$ be the exterior derivative of $f$ along  the leaves. Applying Stoke's theorem to the differential form $\frac{1}{2\pi i} f^{-1} d_{\Fh} f \wedge \omega_{\phi}$ on the complement in $\oX$ of disjoint open tubular neighborhoods of the periodic orbits and of $\varepsilon$-neighborhoods $U_{\varepsilon} (K)$ of the compact leaves $K$, Kopei obtained the following formula for small $\varepsilon > 0$ in \cite{Ko1}
\begin{equation}
\label{eq:7}
\sum_{\gamma} \ord_{\gamma} f \; l (\gamma) - \sum_x \frac{1}{2 \pi i} \int_{\partial U_{\varepsilon} (K_x)} f^{-1} d_{\Fh} f \wedge \omega_{\phi} = 0 \; .
\end{equation}
Here $\gamma$ runs over the periodic orbits and $\ord_{\gamma} f \in \Z$ is defined to be the order of the meromorphic function $f \, |_{F}$ in a point $z \in F \cap \gamma \neq \emptyset$. Since $\ord_z (f \, |_F) \in \Z$ varies continuously with $z \in \gamma$, this is independent of the leaf $F$ and the chosen point $z \in F \cap \gamma$. Moreover $x$ runs over the fixed points and $K_x$ is the compact leaf containing $x$. Comparing \eqref{eq:5} and \eqref{eq:7}, we see that in the analogy prime numbers $p$ correspond to periodic orbits $\gamma$ where $\log p$ corresponds to $l(\gamma)$. Incidentally, note that if we had a bijection $p \leftrightarrow \gamma$ with $\log p = l (\gamma)$, then the Riemann zeta function would be a Ruelle zeta function
\[
\zeta (s) = \prod_{\gamma} (1 - e^{-s l (\gamma)})^{-1} \quad \text{for} \; \RRe s > 1 \; .
\]
The question whether a natural dynamical system with this property exists is quite old and we will discuss our progress in this direction in the last section. Back to our analogy. From \eqref{eq:5} and \eqref{eq:7} we also see that the archimedean absolute value $|\;|_{\infty}$ of $\overline{\spec \Z}$ should correspond to a fixed point $x$. More generally, there is an analogue of \eqref{eq:5} for number fields $K$. It shows that the prime ideals $0 \neq \ep \in \spec \eo_K$ correspond to periodic orbits $\gamma$ where $\log N \ep \ent l (\gamma)$ and the finitely many archimedean absolute values correspond to the finitely many fixed points. Note that every periodic orbit $\gamma$ comes with embeddings
\[
\R / l (\gamma) \Z \hookrightarrow X \; \, \; t \mod l (\gamma) \mapsto \phi^t (x)
\]
for the choices of points $x \in \gamma$. The embeddings for different choices of $x$ are isotopic (via $\phi^t$ for a suitable $t$) and therefore define the same knot in $\oX$. Thus the dynamical systems analogy between primes and periodic orbits refines the previous analogy between primes and knots.

There are several further analogies between $\overline{\spec \Z}$ and systems of the form $(\oX , \Fh , \phi^t)$. For example, Hilbert reciprocity has a dynamical analogue as shown in \cite{MKNT}. Lichtenbaum's conjecture on the order of vanishing and the leading coefficient of the Hasse-Weil zeta function of a regular algebraic scheme over $\spec \Z$ in terms of Weil-\'etale cohomology can be translated almost verbally to the dynamical context, where it can be proved under certain conditions, c.f. \cite{D2}, also \cite{D3}. The proof uses work of \'Alvarez-Lopez and Kordyukov and the Cheeger-M\"uller theorem. 

Here is another analogy which was found earlier \cite[3.5 Corollary]{D4}. The ``explicit formulas of analytic number theory'' are an equality of two distributions on the real line, c.f. \cite{W}. Restricted to $\R^{> 0}$ they have a particularly simple form: In the space of distributions $\Dh' (\R^{> 0})$ we have
\begin{equation}
\label{eq:8} 
1 - \sum_{\rho} e^{t \rho} + e^t = \sum_p \log p \sum_{k \ge 1} \delta_{k \log p} + (1 - e^{-2t})^{-1} \; .
\end{equation}
Here $\rho$ runs over the zeroes of $\zeta (s)$ in $0 < \RRe s < 1$, for $\alpha \in \C$ the locally integrable function $e^{t \alpha}$ is viewed as a distribution, and $\delta_x$ is the Dirac delta distribution supported in $x$. Evaluating on a test function $\varphi \in C^{\infty}_c (\R^{> 0})$ and setting $\Phi (\alpha) = \int_{\R} e^{t \alpha} \varphi (t) \, dt$ formula \eqref{eq:8} takes the more familiar form
\[
\Phi (0) - \sum_{\rho} \Phi (p) + \Phi (1) = \sum_p \log p \sum_{k \ge 1} \varphi (k \log p) + \int^{\infty}_0 \frac{\varphi (t)}{1 - e^{-2t}} \, dt \; .
\]
For the foliation analogue, assume that there exists a metric $g_{\Fh}$ of $T\Fh$ such that $\phi^t$ acts with conformal factor $e^{\alpha t}$ for some $\alpha \in \R$. Then the spectrum of the infinitesimal generator $\theta$ of the group of operators $(\phi^t)^*$ for $t \in \R$ on the space of global leafwise harmonic $L^2$-forms $\ker \Delta^1_{\Fh , (2)}$ consists of eigenvalues $\rho$. If $\dim \oX = 3$ and there are not compact leaves, so that $\oX = X$ and the flow is everywhere transverse to the leaves, then the following formula holds in $\Dh' (\R^{> 0})$ for certain signs $\pm$ which can easily be made explicit
\begin{equation}
\label{eq:9}
1 - \sum_{\rho} e^{t\rho} + e^{\alpha t} = \sum_{\gamma} l (\gamma) \sum_{k \ge 1} \pm \delta_{kl (\gamma)} \; .
\end{equation}
Moreover we have $\RRe \rho = \frac{\alpha}{2}$ for all the eigenvalues $\rho$. Equation \eqref{eq:9} is a transversal index theorem for the $\R$-action $\phi^t$ and the complex $(\Lambda^{\hullet} T^* \Fh , d_{\Fh})$ which is elliptic in the leaf direction. For fixed $t$ the operator $\phi^{t*}$ on $\ker \Delta^i_{\Fh , (2)}$ is not trace class if $\ker \Delta^i_{\Fh , (2)}$ is infinite dimensional. However for $\varphi \in C^{\infty}_c (\R)$ the trace of the mollified operator $\int_{\R} \phi^{t*} \varphi (t) \, dt$ on $\ker \Delta^i_{\Fh , (2)}$ exists and by our assumptions we have
\begin{align*}
\langle \Tr (\phi^* \mid \ker \Delta^i_{\Fh, (2)} , \varphi \rangle & := \Tr \int_{\R} \phi^{t*} \varphi (t) \, dt \\
& = \sum_{\rho} \int_{\R} e^{t \theta} \varphi (t) \, dt \\
& = \langle \sum_{\rho} e^{t\rho} , \varphi \rangle \;.
\end{align*}
The only eigenvalues $\rho$ of $\theta$ on $\ker \Delta^i_{\Fh , (2)}$ for $i = 0,2$ are $0$ and $\alpha$ by our conditions. Thus the left hand side of \eqref{eq:9} may be viewed as the transversal index defined by Atiyah for compact Lie group actions and by H\"ormander in general. It is given by the following Euler-characteristic
\[
\sum_i (-1)^i \Tr (\phi^* \mid \ker \Delta^i_{\Fh, (2)}) \in \Dh' (\R) \; .
\]
Incidentally, $\ker \Delta^i_{\Fh, (2)}$ may also be interpreted as the maximal Hausdorff quotient of the leafwise $L^2$-cohomology
\[
H^i_{\Fh , (2)} (X) = \ker d^i_{\Fh , (2)} / \imm d^{i-1}_{\Fh , (2)} \; .
\]
Formula \eqref{eq:9} does not contain a term corresponding to $(1 - e^{-2 t})^{-1}$ in \eqref{eq:8} because we assumed that $\phi^t$ had no fixed points. If we allow fixed points, then the distributional trace defined above may no longer exist. An extension of transverse index theory to such a more general situation has only recently been accomplished by \'Alvarez-Lopez, Leichtnam and Kordyukov, c.f. \cite{ALK2}, \cite{ALK3}. It is much more involved then the transversal case and leads to formulas which are quite similar to the explicit formulas of number theory, even as distributions on all of $\R$.

One can show that the conformal factor $e^{\alpha t}$ for a metric $g_{\Fh}$ as above necessarily has to be $1$, i.e. $\alpha = 0$ which implies that the eigenvalues $\rho$ of $\theta$ on $\ker \Delta^1_{\Fh , (2)}$ have real part $\RRe \rho = 0$. By comparison with the Riemann hypothesis, we would obviously want $\alpha = 1$ to be a possibility in the geometric analogue $(\oX , \Fh , \phi^t)$. This is only one instance which shows that the class of smooth compact manifolds is too restrictive to actually rewrite the explicit formulas of analytic number theory as a transversal index theorem. One can show that $\alpha = 1$ can be achieved if for $\oX$ we allow the local structure (totally disconnected) $\times$ ($3$-dimensional ball), c.f. \cite{L}. In \cite{G} and \cite{MS} it is explained how to do analyis of PDE on such spaces. Analogies of Arakelov theory with foliated dynamical systems were obtained in \cite{Ko2}. Recently we found an argument why a possible complex valued Weil-type cohomology theory for arithmetic curves $\overline{\spec \eo_K}$ cannot have a functorial real structure, \cite{D5}. Since leafwise cohomology groups or the spaces of global leafwise harmonic forms always have natural real structures, this shows that our analogies are not perfect in this regard. We are missing a fundamental ``twist'' excluding real structures on cohomology somewhere. 
\section{Dynamical systems for arithmetic schemes} \label{sec:4}

In this section we recall the construction of ``foliated'' topological dynamical systems for arithmetic schemes given in \cite{D6} and explain their basic properties. We first found an \textit{extrinsic} construction of the typical leaf with its Poincar\'e return maps guided by the idea to use Frobenius elements in Galois groups to generate periodic orbits. The following more general and conceptual \textit{intrinsic} definition using $\C$-valued points of rational Witt spaces came later after understanding the work of Kucharczyk and Scholze \cite{KSc}. The previous extrinsic definition then became a theorem, namely Theorem \ref{t41} below. Our dynamical systems have some but by no means all of the properties that we expect from the analogies in section \ref{sec:3}. We view them as a first but important step in our quest to apply analysis on dynamical systems to number theory. Since arithmetic over $p$-adic fields is much better understood than over number fields we also studied a $p$-adic version of our construction. In that situation there is a very natural modification of the dynamical system and it turned out to be closely related to the Fargues-Fontaine curve, one of the fundamental objects in $p$-adic Hodge theory.

Before we can define rational Witt spaces, we have to recall the rational Witt vectors $W_{\rat} (R)$ of a commutative unital ring $R$. As a set, $W_{\rat} (R)$ consists of the rational functions $f$ among the power series $f \in R [[t]]$, with $f (0) = 1$, i.e. $f = P / Q$ with $P, Q \in R [t] , P (0) = 1 = Q (0)$. Addition in $W_{\rat} (R)$ is defined to be multiplication of rational functions. Following Almkvist \cite{Almkvist} we write $f \in W_{\rat} (R)$ as a quotient
\[
f = \frac{\det (1 - t \varphi \mid V)}{\det (1 - t \psi \mid W)} \; ,
\]
where $V$ and $W$ are projective $R$-modules of finite rank equipped with endomorphisms $\varphi$ and $\psi$. If $\tf$ has a corresponding expression then the sum $f + \tf$ in $W_{\rat} (R)$ will be represented by $(\varphi \oplus \tvarphi , V \oplus \tV)$ and $(\psi \oplus \tpsi , W \oplus \tW)$. By definition the product of $f$ and $\tf$ in $W_{\rat} (R)$ is represented by $(\varphi \otimes \tvarphi , V \otimes \tV)$ and $(\psi \otimes \tpsi , W \otimes \tW)$. One can check that this is independent of all choices, and in fact the product in $W_{\rat} (R)$ coincides with the product in the big Witt ring $W (R) = 1 + tR [[t]]$. This leads to the $K$-theoretical description of $W_{\rat} (R)$ in \cite{Almkvist}. Now let $X$ be a scheme, for example $X = \spec R$ and let $W_{\rat} (\Oh_X)$ be the sheafification of the Zariski presheaf $U \mapsto W_{\rat} (\Oh_X (U))$. For all $x \in X$ we have $W_{\rat} (\Oh_X)_x = W_{\rat} (\Oh_{X,x})$. We call the ringed space
\[
W_{\rat} (X) = (X_{\top} , W_{\rat} (\Oh_X)) \; ,
\]
the rational Witt space of $X$, where $X_{\top}$ is the underlying topological space of $X$. If $S$ is another scheme we define a morphism $P : S \to W_{\rat} (X)$ as a morphism of ringed spaces $(P , P^{\sharp})$ such that for all $s \in S$ there is a (uniquely determined) ring homomorphism $\tP^{\sharp}_s$ making the following diagram commutative
\begin{equation}
\label{eq:10}
\xymatrix{W_{\rat} (\Oh_{X , f (s)}) \ar[r]^-{P^{\sharp}_s} \ar@{>>}[d] & \Oh_{S,s} \ar@{>>}[d] \\
W_{\rat} (\kappa (f (s))) \ar[r]^-{\tP^{\sharp}_s} & \kappa (s) \; .
}
\end{equation}
Here $\kappa$ denotes the residue field. The stalks $W_{\rat} (\Oh_X)_s$ are not local rings in general and $W_{\rat} (X)$ is not a locally ringed space. Condition \eqref{eq:10} replaces the usual locality condition in scheme theory. The set of $S$-valued points of $W_{\rat} (X)$ is
\[
W_{\rat} (X) (S) = \Mor (S , W_{\rat} (X)) \; .
\]
If $S = \spec A$ is affine we set $W_{\rat} (X) (A) = W_{\rat} (X) (S)$. The canonical ring homomorphisms
\[
W_{\rat} (R) \longrightarrow R \; , \; f \longmapsto -f' (0) / f (0)
\]
induce a canonical morphism
\[
X \longrightarrow W_{\rat} (X) \; .
\]
By composition we get an injection
\[
X (S) \hookrightarrow W_{\rat} (X) (S) \; .
\]
Hence $W_{\rat} (X) (S)$ contains the classical $S$-valued points of $X$. 

Morphisms $W_{\rat} (Y) \to W_{\rat} (X)$ between two rational Witt spaces are defined similarly as above. Every morphism $\alpha : Y \to X$ of schemes induces a morphism $W_{\rat} (\alpha) : W_{\rat} (Y) \to W_{\rat} (X)$ and one obtains a faithful but not fully faithful functor from schemes to rational Witt spaces. For example, the commuting Frobenius endomorphisms $F_{\nu}$ for $\nu \ge 1$ of Witt vector theory induce commuting Frobenius endomorphisms $F_{\nu}$ on $W_{\rat} (X)$. We have $F_{p ^n} = W_{\rat} (F^n)$ if $X$ is an $\F_p$-scheme and $F : X \to X$ is the absolute Frobenius, but in general the $F_{\nu}$ are not induced by morphisms of schemes. The multiplicative monoid $\Nh$ of positive integers therefore acts on $W_{\rat} (X)$ and hence also on $W_{\rat} (X) (S)$ for all $S$. We now describe how global sections $f$ of $X$ give rise to $A$-valued functions on $W_{\rat} (X) (A)$ using the multiplicative (Teichmüller-)map
\[
[\;] : R \longrightarrow W_{\rat} (R) \quad \text{with} \; [r] = 1 - rt \; .
\]
Namely, for $f \in \Gamma (X , \Oh_X)$ consider its images $[f]$ in $W_{\rat} (\Oh_X (X))$ and hence in $W_{\rat} (\Oh_X) (X)$. A point $P \in W_{\rat} (X) (\C)$ gives a morphism of sheaves $P^{\sharp} : W_{\rat} (\Oh_X) \to P_* \Oh_{\spec \C}$. We define the value of $f$ in $P$ to be the image of $[f] \in W_{\rat} (\Oh_X) (X)$ in $(P_* \Oh_{\spec \C}) (X) = \C$ under this map. In this way we obtain a multiplicative map from $\Gamma (X , \Oh_X)$ into the $\C$-algebra of $\C$-valued function on $W_{\rat} (X) (\C)$.

For normal schemes $X$ with countable function field we defined a natural topology on $W_{\rat} (X) (\C)$ in \cite[section 7]{D6} and the previous construction gives a multiplicative map
\begin{equation}
\label{eq:11}
[\;] : \Gamma (X , \Oh) \longrightarrow C^0 (W_{\rat} (X) (\C) , \C) \; .
\end{equation}
For $X = \spec \Z$ this map interprets numbers i.e. elements of $\Gamma (\spec \Z , \Oh) = \Z$ as highly non-trivial $\C$-valued continuous functions on the infinite dimensional connected Hausdorff space $W_{\rat} (\spec \Z) (\C)$. The following result which follows from \cite[Lemma 4.9]{KSc} makes the points of $W_{\rat} (X) (\C)$ somewhat more explicit but even for $X = \spec \Z$ we do not have a satisfactory understanding of that space. We change the notation somewhat.

\begin{theorem}
\label{t41}
Let $\eX_0$ be a normal scheme with function field $K_0$ and let $\eX$ be the normalization of $\eX_0$ in an algebraic closure $K$ of $K_0$. Let $G = \Aut (K / K_0)$ be the absolute Galois group of $K_0$ and set
\[
\deX (\C) = \{ (\ex , \oP^{\times}) \mid \ex \in \eX , \oP^{\times} \in \Hom (\kappa (\ex)^{\times} , \C^{\times}) \} \; .
\]
Then there is a natural $\Nh$-equivariant bijection
\[
W_{\rat} (\eX_0) (\C) \silo \deX_0 (\C) := \deX (\C) / G \; .
\]
Here $\sigma \in G$ acts on $\deX (\C)$ via $(\ex , \oP^{\times})^{\sigma} = (\ex^{\sigma} , \oP^{\times} \verk \sigma)$ and $\nu \in \Nh$ acts $G$-equivariantly via $F_{\nu} (\ex , \oP^{\times}) = (\ex , \oP^{\times} \verk (\;)^{\nu})$. 
\end{theorem}

The space $W_{\rat} (\eX_0) (\C) = \deX_0 (\C)$ is almost the typical leaf of the ``foliated'' dynamical system that we will construct. The monoid $\Nh$ acts on $\deX_0 (\C)$ by injective maps and we need to invert them to get a group action by $(\Q^{> 0} , \cdot)$ via homeomorphisms. This is done by forming the colimit space over $\Nh$ viewed as a poset ordered by divisibility
\[
\ceX_0 (\C) = \colim_{\Nh} \deX_0 (\C) \; .
\]
Note that on the level of functions we have a limit over Frobenii, which reminds of the tilting process in $p$-adic geometry. Set
\[
X_0 = (\ceX_0 (\C) \times \R^{> 0}) / \Q^{> 0}
\]
where $\Q^{> 0}$ acts diagonally. Let $t \in \R$ act on $X_0$ by setting $\phi^t [P , u] = [P , e^t u]$. The $1$-codimensional ``foliation'' $\Fh$ has leaves the images of $\ceX_0 (\C) \times \{ u \}$ in $X_0$ for $u \in \R^{> 0}$. It is everywhere transversal to the flow and each $\phi^t$ maps leaves to leaves. In general, the dynamical system $(X_0 , \phi^t)$ has too many periodic orbits, since the $\Nh$-space $W_{\rat} (\eX_0) (\C)$ does not know enough about the addition in $\Oh_{\eX_0}$. In the local $p$-adic situation below, we know the right modification to make. However in the global case presently we can only impose an ``admissible'' condition $\Eh$ on the characters $\oP^{\times} : \kappa (\ex)^{\times} \to \C^{\times}$ in the description of Theorem \ref{t41}. What we know for certain is that the restrictions of $\oP^{\times}$ to $\mu (\kappa (\ex))$ must have finite kernels (condition $\Eh_{\tors})$ and that for $\car \kappa (\ex) > 0$ the image of $\oP^{\times}$ must not be torsion unless $\kappa (\ex)^{\times}$ is itself torsion. For example $\Eh$ can be the conditions that $\ker \oP^{\times}$ is always finite resp. finitely generated. We refer to \cite[section 4]{D6} for a detailed discussion. With the obvious modifictions we get a $G \times \Nh$-space $\deX (\C)_{\Eh}$, an $\Nh$-space $\deX_0 (\C)_{\Eh}$ and a ``foliated'' dynamical system $(X_{0 \Eh} , \Fh , \phi^t)$. 

\begin{theorem}
\label{t42}
Let $\eX_0$ be normal of finite type over $\spec \Z$, e.g. $\eX_0 = \spec R_0$ for an integrally closed finitely generated ring $R_0$. Then we have
\[
\{ x_0 \in X_{0\Eh} \mid \phi^t (x_0) = x_0 \quad \text{for some} \; t > 0 \} = \coprod_{\ex_0} \Gamma_{\ex_0} \; .
\]
Here $\ex_0$ runs over the closed points of $\eX_0$, i.e. the maximal ideals $\emm_0$ of $R_0$ if $\eX_0 = \spec R_0$. The compact subsets $\Gamma_{\ex_0} \subset X_0$ consist of periodic orbits of length $\log N \ex_0$ where $N \ex_0 = |\kappa (\ex_0)|$ ($= |R_0 / \emm_0|$) and they are pairwise disjoint. In fact $\Gamma_{\ex_0}$ is a fibre space over the compact group $\Aut (\oF^{\times}_p) / \Aut (\oF_p)$ where $p = \car \kappa (\ex)$ with fibres the compact orbits in $\Gamma_{\ex_0}$. 
\end{theorem}

The theorem asserts that the closed points of $\eX_0$ e.g. the prime numbers $p$ if $\eX_0 = \spec \Z$ correspond not to individual periodic orbits $\gamma$ as in the analogies of section \ref{sec:3} but to compact packets of periodic orbits all of which have length $\log N \ex_0$ i.e. $\log p$ if $\eX_0 = \spec \Z$. This is reminiscient of the invariant tori of Hamiltonian dynamics.

It follows from the existence of Frobenius elements in the Galois group $G$ that closed points give periodic orbits in $X_{0 \Eh}$. It is more difficult to show that any periodic orbit comes from a closed point i.e. lies in $\Gamma_{\ex_0}$ for some $\ex_0$, c.f. \cite[Theorems 5.2 and 6.1]{D6}.

The proof of the following result requires de Jong's theory of alterations if $\dim \eX_0 \ge 2$ and some approximation theorems from number theory.

\begin{theorem}
\label{t43}
Let $\eX_0$ be an integral normal scheme which is flat of finite type over $\spec \Z$. Then the spaces $X_0$ and $X_{0 \Eh}$ are connected. In fact they are almost pathwise connected: for any two points $x_0$ and $x'_0$ in $X_{0\Eh}$ and any two neighborhoods $U_0 \ni x_0$ and $U'_0 \ni x'_0$ there are points $y_0 \in U_0$ and $y'_0 \in U'_0$ which can be connected by a continuous path.
\end{theorem}

Another result, which is stronger than connectedness asserts that the leafwise cohomology group $H^0_{\Fh} (X_{0 \Eh})$ is one dimensional. It is essentially the group of global sections of the sheaf of continuous $\R$-valued functions on $X_{0\Eh}$ which are locally constant along the leaves of $\Fh$, c.f. \cite[section 10]{D6}. 

The space $X_{0 \Eh}$ is infinite dimensional if $\dim \eX_0 \ge 1$ and one could hope that the sub-dynamical system obtained as the closure of the union of all its compact orbits might be significantly smaller. However, this is not the case as follows from \cite[Theorem 8.2]{D6}. The closure is the subsystem obtained by replacing $\deX (\C)_{\Eh}$ in the previous constructions with the subspace of pairs $(\ex , \oP^{\times})$ with $\oP^{\times} : \kappa (\ex)^{\times} \to S^1$ a unitary character. For $\dim \eX_0 \ge 2$ this is conditional on a result in Diophantine approximation which should be provable and which is known for $\dim \eX_0 = 1$ and $\eX_0$ flat over $\spec \Z$.

The following result makes the structure of the dynamical system $X_0$ clearer. Consider the natural projection $\deX (\C) \to \eX$ mapping $(\ex , \oP^{\times})$ to $\ex$. Let $\deX (\C)_p$ resp. $\deX (\C)_{\Q}$ be the fibres of the composition $\deX (\C) \to \eX \to \spec \Z$ over $(p)$ resp. $(0)$. Consider the $G$-invariant subspace
\[
\deX (\C)_{\inn} = \{ (\ex , \oP^{\times}) \in \deX (\C) \mid \oP^{\times} \, |_{\mu (\kappa (\ex))} \; \text{is injective} \, \}
\]
and set
\[
\deX (\C)_{p,\inn} = \deX (\C)_{\inn} \cap \deX (\C)_p \quad \text{and} \quad \deX (\C)_{\Q , \inn} = \deX (\C)_{\inn} \cap \deX (\C)_{\Q} \; .
\]
These are subspaces of $\deX (\C)$ and the quotient topologies on 
\[
\deX_0 (\C)_{p , \inn} = \deX (\C)_{p , \inn} / G \quad \text{and} \quad \deX_0 (\C)_{\Q , \inn} = \deX (\C)_{\Q , \inn} / G
\]
agree with the subspace topologies within $\deX_0 (\C)$. For any $p$, the Frobenius endomorphism $F_p$ of $\deX_0 (\C)$ restricts to a homeomorphism of $\deX_0 (\C)_{p , \inn}$. Hence $p^{\Z} \subset \Q^{> 0}$ acts on $\deX_0 (\C)_{p , \inn}$ via $p \leftrightarrow F_p$.

\begin{theorem}
\label{t44}
The following canonical $\phi^t$-equivariant map is a continuous bijection (and similarly after imposing an admissible condition $\Eh$)
\[
\deX_0 (\C)_{\Q , \inn} \times \R^{> 0} \amalg \coprod_p \deX_0 (\C)_{p , \inn} \times_{p^{\Z}} \R^{> 0} \silo X_0 = \ceX_0 (\C)_{\Eh_{\tors}} \times_{\Q^{> 0}} \R^{> 0} \; .
\]
\end{theorem}

In view of Theorem \ref{t43} the map in the theorem is not a homeomorphism since $X_0$ is connected, whereas the left hand side is disconnected. 

We now discuss the relation of our construction with the work of Kucharczyk and Scholze \cite{KSc}. Consider a field of characteristic zero containing all roots of unity, and fix an injective homomorphism $\iota : \mu (F) = \mu (\onF) \hookrightarrow \C^{\times}$. The multiplicative Teichm\"uller map $[\;] : F^{\times} \hookrightarrow W_{\rat} (F)$ and the multiplicative map $\iota$ give ring homomorphisms
\[
\Z [\mu (F)] \longrightarrow W_{\rat} (F) \quad \text{and} \quad \Z [\mu (F)] \longrightarrow \C
\]
on the group ring of $\mu (F)$. Set
\[
\eX_F = \spec (W_{\rat} (F) \otimes_{\Z [\mu (F)]} \C) \; .
\]
One can show that the connected components of $W_{\rat} (\spec F) (\C)$ are parametrized by the embeddings $\mu (F) \hookrightarrow \C^{\times}$ and that $\eX_F (\C)$ is the component corresponding to $\iota$. For arbitrary connected pointed topological spaces $(Z, z)$ Kucharczyk and  Scholze defined a pro-finite fundamental group $\pi^{\et}_1 (Z , z)$ which classifies the finite coverings of $Z$. One of their main results is an isomorphism of $\pi^{\et}_1 (\eX_F (\C) , x)$ with the absolute Galois group of $F$. They also calculate the sheaf cohomology of $\eX_F (\C)$ with coefficients in $\Z / n$ and $\Q$. For $\Z / n$ they prove that it is isomorphic to the Galois cohomology of $F$ with $\Z / n$-coefficients using the work of Rost and Voevodsky on the Bloch-Kato conjecture. For $\Q$-coefficients the cohomology of $\eX_F (\C)$ is the Galois fixed part of $\Lambda^{\hullet} (\onF^{\times} \otimes \Q)$. A more natural group would be the Milnor $K$-group of $F$ tensored with $\Q$. The absence of the Steinberg relations in rational cohomology is an indication that the space $\eX_F (\C)$ and hence also our space $W_{\rat} (X) (\C)$ do not encode enough information about the additive structure of $F$ resp. $\Oh_X$. 

In \cite{DW} we introduced and studied a pro-algebraic fundamental group $\pi_E (Z, z)$ over any field $E$ for pointed connected topological spaces $(Z,z)$. Its maximal pro-\'etale quotient is $\pi^{\et}_1 (Z,z)$ viewed as a group scheme over $E$. For the ``right'' version of $\eX_F (\C)$ discussed in the introduction of \cite{KSc}, we expect the pro-algebraic fundamental group to be deeply related to the motivic Galois group for motives over $F$ with coefficients in $E$. In \cite[section 5]{DW} we calculated $\pi_E (\eX_F (\C) , x)$ for algebraically closed fields $E$ of characteristic zero. Its connected component is commutative and its unipotent part is related to extensions of Tate motives. For the reductive part we could see no relation with motives, i.e. to the Serre group. In our opinion this just shows that $\eX_F (\C)$ and our spaces $W_{\rat} (\eX_0) (\C)$ are only first steps (but in the right direction). 

For inspiration, we looked at the much better understood situation where $\eX_0$ is a normal scheme of finite type over $\spec \Z_p$ instead of $\spec \Z$. In order to compare with $p$-adic Hodge theory it is reasonable to study $W_{\rat} (\eX_0) (\eo)$ where $\eo$ is the valuation ring in the completion $\C_p$ of $\oQ_p$. We proved that
\[
W_{\rat} (\eX_0) (\eo) = W_{\rat} (\eX) (\eo)/ G \; .
\]
Here, as before $\eX$ is the normalization of $\eX_0$ in $\oK_0$, with $K_0$ the function field of $\eX_0$ and $G$ its absolute Galois group. The points of $W_{\rat} (\eX) (\eo)$ are certain diagrams determined by points $\ex , \ey \in \eX$ with $\ey \in \overline{\{ \ex \}}$ and multiplicative maps $P_{\ey} : \Oh_{\overline{\{ \ex \}} , \ey} \to \eo$ sending $1$ to $1$ and $0$ to $0$. Here we know how to modify our constructions in order to get the right $F_p$-dynamical system: We have to impose the condition that the $\mod p$ reduction of the multiplicative map $P_{\ey}$ i.e. the composition 
\begin{equation}
\label{eq:12}
\Oh_{\overline{\{ \ex \}} , \ey} \xrightarrow{P_{\ey}} \eo \longrightarrow \eo / p
\end{equation}
is \textit{additive}. The technical details are more involved but this is the heart of the matter. It now turns out that the elements of Fontaine's $A_{\inf}$-rings become ordinary $\eo$-valued functions on the resulting sub-dynamical systems. Judith Lutz has proved that this representation of $A_{\inf}$ by functions is faithful. Also, e.g. for $\eX_0 = \spec \Z_p$ we obtained a natural bijection of the subsystem with the points of the Fontaine-Fargue curve. This uses the work of Fontaine-Wintenberger in the reformulation by Scholze and Kedlaya. One can reformulate additivity of \eqref{eq:12} in terms of absolute values but it is still unclear what the right analogue is to impose on $W_{\rat} (\eX) (\C)$ for $\eX_0 / \Z$. 

\end{document}